\documentclass[11pt]{article}
\usepackage{mathptmx}
\usepackage{times}
\usepackage{latexsym,epsfig,bm}
\usepackage{amssymb}
\usepackage{color}
\usepackage{amsthm}
\usepackage{mathrsfs}

\def\version{April 27, 2007}

\DeclareSymbolFont{EUR}{U}{eur}{m}{n}
\SetSymbolFont{EUR}{bold}{U}{eur}{b}{n}
\DeclareSymbolFontAlphabet{\eur}{EUR}

\DeclareSymbolFont{EUB}{U}{eur}{b}{n}
\SetSymbolFont{EUB}{bold}{U}{eur}{b}{n}
\DeclareSymbolFontAlphabet{\eub}{EUB}

\DeclareSymbolFont{AMSb}{U}{msb}{m}{n}
\DeclareSymbolFontAlphabet{\mathbb}{AMSb}

\newcommand{\beqn}{\begin{eqnarray}}
\newcommand{\eeqn}{\end{eqnarray}}
\newcommand{\be}{\begin{equation}}
\newcommand{\ee}{\end{equation}}

\newcommand{\na}{\nabla}

\newcommand\supp{\mathop{\rm supp}}

\newcommand{\at}[1]{\vert\sb{\sb{#1}}}

\def\R{{\mathbb R}}

\newcommand\C{{\mathbb C}}

\newcommand{\Abs}[1]{\left\vert#1\right\vert}
\newcommand{\abs}[1]{\vert #1 \vert}
\newcommand{\Norm}[1]{\left\Vert #1 \right\Vert}
\newcommand{\norm}[1]{\Vert #1 \Vert}

\DeclareMathSymbol{\varGamma}{\mathord}{letters}{"00}
\DeclareMathSymbol{\varDelta}{\mathord}{letters}{"01}
\DeclareMathSymbol{\varTheta}{\mathord}{letters}{"02}
\DeclareMathSymbol{\varLambda}{\mathord}{letters}{"03}
\DeclareMathSymbol{\varXi}{\mathord}{letters}{"04}
\DeclareMathSymbol{\varPi}{\mathord}{letters}{"05}
\DeclareMathSymbol{\varSigma}{\mathord}{letters}{"06}
\DeclareMathSymbol{\varUpsilon}{\mathord}{letters}{"07}
\DeclareMathSymbol{\varPhi}{\mathord}{letters}{"08}
\DeclareMathSymbol{\varPsi}{\mathord}{letters}{"09}
\DeclareMathSymbol{\varOmega}{\mathord}{letters}{"0A}

\theoremstyle{plain}
\newtheorem{theorem}{Theorem}[section]

\newtheorem{lemma}[theorem]{Lemma}
\newtheorem{proposition}[theorem]{Proposition}
\newtheorem{corollary}[theorem]{Corollary}

\theoremstyle{definition}

\theoremstyle{remark}

\newcommand\const{\mathop{{\rm const}}}

\makeatletter\@addtoreset{equation}{section}
\makeatother

\renewcommand{\Re}{\mathop{\eur{R\hskip -1pt e}}\nolimits}
\renewcommand{\Im}{\mathop{\eur{I\hskip -1pt m}}\nolimits}

\usepackage{epsfig}
\usepackage{graphicx}

\newenvironment{equations}{\begin{eqnarray}}{\end{eqnarray}}

\begin{document}

\title{Global well-posedness
for the Schr\"odinger equation
\\
coupled to a nonlinear oscillator
}

\author{
{\sc Alexander Komech}
\footnote{
On leave from Department of Mechanics and Mathematics,
Moscow State University, Moscow 119992, Russia.
Supported in part
by  DFG grant 436\,RUS\,113/615/0-1, FWF grant P19138-N13,
Max-Planck Institute for Mathematics in the Sciences (Leipzig)
and Alexander von Humboldt Research Award (2006).}
\\
{\it\small Faculty of Mathematics, Vienna University, Wien A-1090, Austria}
\\
{\sc Andrew Komech}
\footnote{
Supported in part
by Max-Planck Institute for Mathematics in the Sciences (Leipzig)
and by the NSF Grants DMS-0434698 and DMS-0600863.
}
\\
{\it\small
Mathematics Department, Texas A\&M University,
College Station, TX, USA}
}
\date{\version}

\maketitle

\begin{abstract}
The Schr\"odinger equation
with the nonlinearity concentrated at a single point
proves to be an interesting and important model
for the analysis of long-time behavior of solutions,
such as the asymptotic stability of solitary waves
and properties of weak global attractors.
In this note, we prove global well-posedness
of this system in the energy space $\scriptstyle{H\sp 1}$.
\end{abstract}


\section{Introduction and main results}
\label{sect-introduction}

We are going to prove the well-posedness in
$H\sp 1$
for the nonlinear Schr\"odinger equation
with the nonlinearity concentrated at a single point:
\begin{equation}\label{snls}
i\dot\psi(x,t)
=-\psi''(x,t)-\delta(x)F(\psi(0,t)),
\qquad
\quad x\in\R,
\end{equation}
where the dots and the primes
stand for the partial derivatives in $t$ and $x$,
respectively.
The equation describes the Schr\"odinger field
coupled to a nonlinear oscillator.
This equation is a convenient playground for developing the tools
for the analysis of long-time behavior of solutions to
$U(1)$-invariant Hamiltonian systems with dispersion.
The asymptotic stability of the solitary manifold
for equation (\ref{snls})
has been considered in \cite{bkks}.
Here we complete this result,
giving the proof of the global well-posedness of (\ref{snls})
in the energy space.

Let us mention that for the Klein-Gordon equation
with the nonlinearity of the same type
the global attraction was addressed in 
\cite{ubk-cr},
\cite{ubk-arma}.

We assume that
\begin{equation}\label{FU}
F(\psi)=-\nabla\sb\psi U(\psi),
\qquad\psi\in\C,
\end{equation}
for some real-valued potential $U\in C\sp 2(\C)$,
where $\nabla\sb\psi$ is the real derivative
with respect to $(\Re\psi,\Im\psi)$.
Equation (\ref{snls}) is a Hamiltonian system
with the Hamiltonian
\begin{equation}
\mathscr{H}(\psi)
=\int\sb{\R}
\frac{\abs{\psi'(x)}^2}{2}
\,dx
+U(\psi(0)),
\qquad
\psi\in H^1=H\sp 1(\R).
\end{equation}
The Hamiltonian form of (\ref{snls}) is
\begin{equation}\label{Ham}
\dot\Psi=J D \mathscr{H}(\Psi),
\end{equation}
where
\begin{equation}
\Psi=\left[\begin{array}{c}\Re\psi\\\Im\psi\end{array}\right],
\qquad
J=\left[\begin{array}{cc}0&1\\-1&0\end{array}\right],
\end{equation}
and $D \mathscr{H}$ is the Fr\'echet derivative on the Hilbert space $H^1$.
The value of the Hamiltonian functional
is conserved for classical finite energy solutions of (\ref{snls}).
We assume that equation (\ref{snls}) possesses
$U(1)$-symmetry,
thus requiring that
\begin{equation}\label{Uu}
U(\psi)=u(\abs{\psi}^2),
\qquad
\psi\in\C.
\end{equation}
It then follows that
$F(0)=0$
and
$F(e^{i s}\psi)=e^{i s}F(\psi)$
for $\psi\in\C$, $s\in\R$,
and that
\begin{equation}\label{def-a}
F(\psi)=a(\abs{\psi}^2)\psi,
\qquad
\psi\in\C,
\qquad
{\rm where}
\quad
a(\cdot)=2 u'(\cdot)\in\R.
\end{equation}
This symmetry implies that
$e^{i\theta}\psi(x,t)$ is a solution to (\ref{snls})
if $\psi(x,t)$ is.
According to the N\"other theorem,
the $U(1)$-invariance leads (formally)
to the conservation of
the charge, given by the functional
\begin{equation}\label{charge}
Q(\psi)=\frac{1}{2}\int\sb{\R}\abs{\psi}^2\,dx.
\end{equation}
We also assume that $U(\psi)$ is such that
\begin{equation}\label{u-is-such}
U(z)\ge A-B\abs{z}^2
\qquad
{\rm with\ some}
\quad
A\in\R,\quad B > 0.
\end{equation}
We will show that
equation (\ref{snls}) is globally well-posed in $H\sp 1$.
We will consider the solutions of class
$\psi\in C\sb b(\R\times\R) $.
All the derivatives in equation (\ref{snls})
are understood in the sense of distributions.

\begin{theorem}[Global well-posedness]
\label{theorem-gwp}
Let the conditions
(\ref{FU}), (\ref{Uu})  and (\ref{u-is-such}) hold with $U\in C^2(\C)$. Then
\begin{enumerate}
\item
For any $\phi\in H^1(\R)$,
the equation for (\ref{snls})
with the initial data $\psi\at{t=0}=\phi$
has a unique solution
$\psi\in C(\R,H^1(\R))$.

\item
The values of the charge and energy functionals
are conserved:
\begin{equation}
\label{star}
Q(\psi(t))=Q(\phi),
\qquad
\mathscr{H}(\psi(t))=\mathscr{H}(\phi),
\qquad t\in\R.
\end{equation}
\item
There exists $\Lambda(\phi)>0$
such that
the following a priori bound holds:
\begin{equation}
\sup\limits\sb{t\in\R}\norm{\psi(t)}\sb{H^1}
\le \Lambda(\phi)<\infty.
\label{a-priori-h1}
\end{equation}
\item
The map
$\eub{U}:\,\psi(0)\mapsto\psi$
is continuous
from $H^1$ to $L\sp\infty([0,T],H^1(\R))
$,
for any $T>0$.
\end{enumerate}
\end{theorem}

\begin{theorem}\label{theorem-holder}
Under conditions of Theorem~\ref{theorem-gwp},
$\psi\in C\sp{(1/4)}(\R\times\R)$.
\end{theorem}


Let us give the outline of the proof.
We need a small preparation first:
We show that, without loss of generality,
it suffices to prove the theorem assuming
that $U$ is uniformly bounded together with its
derivatives.
Indeed, the a priori bounds on the $L\sp\infty$-norm
of $\psi$ imply that the nonlinearity $F(z)$
may be modified for large values of $\abs{z}$.
Then we will prove the existence and
uniqueness of the solution
$\psi\in  C\sb b(\R\times[0,\tau])$,
for some $\tau>0$.
This is accomplished in Section~\ref{sect-cb}.

In Section~\ref{sect-regularized},
we construct
approximate solutions
$\psi\sb\epsilon\in C\sb b(\R,H^1(\R))$
that are solutions to a regularized problem
(the $\delta$-function substituted by its smooth approximations
$\rho\sb\epsilon$, $\epsilon>0$).
On one hand, the approximate solutions have their
energy and charge conserved.
On the other hand,
we will show in Section~\ref{sect-convergence}
that the approximate solutions converge to $\psi(x,t)$
uniformly for $\abs{x}\le R$, $0\le t\le\tau$.

In Section~\ref{sect-wp},
we use the uniform convergence of approximate solutions
to conclude that $\psi\in L\sp\infty([0,\tau],H^1(\R))$
and moreover that $\psi$ could be extended to all $t\ge 0$.
Then we show that
the energy and the charge are conserved.
We will use these conservations
to extend the solution $\psi(x,t)$ for $t\in\R$.
Then we prove that $\psi\in C(\R,H^1(\R))$.

In Section~\ref{sect-holder},
we study the H\"older continuity in time,
showing that $\psi\in C\sp{(1/4)}(\R\times\R)$.

\section{Local well-posedness in $C\sb b$}
\label{sect-cb}

\begin{lemma}\label{lemma-a-priori}
A priori bound
(\ref{a-priori-h1}) follows from (\ref{u-is-such}) and
the energy and charge conservation
(\ref{star}).
\end{lemma}

\begin{proof}
Let $A\in\R$, $B>0$
be constants from (\ref{u-is-such}),
and let $\psi\in H\sp 1(\R)$.
To estimate $\norm{\psi}\sb{H\sp 1}$ in terms of the values
of $Q(\psi)$ and $\mathscr{H}(\psi)$,
we need to control the possibly negative contribution of $U(\psi)$
into $\mathscr{H}(\psi)$.
We achieve this control
by using the inequality
\begin{equation}\label{psi-0}
B\abs{\psi(0)}^2
\le
B\Big[\int\sb\R\hat\psi(k)\frac{dk}{2\pi}\Big]^2
\le
B\int\sb\R \big(B^2+{\textstyle\frac{k^2}{4}}\big)
\abs{\hat\psi(k)}^2\frac{dk}{2\pi}
\cdot
\int\sb R\frac{dk}{2\pi\big(B^2+\frac{k^2}{4}\big)}
=B^2\norm{\psi}\sb{L\sp 2}^2+{\textstyle\frac{1}{4}}\norm{\psi'}\sb{L\sp 2}^2.
\end{equation}
This allows us to write
\begin{equation}
\label{q-h-0}
\mathscr{H}(\psi)
\ge
{\textstyle\frac 1 2}
\norm{\psi'}^2\sb{L\sp 2}
+A-B\abs{\psi(0)}^2
\ge
{\textstyle\frac 1 4}
\norm{\psi'}^2\sb{L\sp 2}
+A-B^2\norm{\psi}\sb{L\sp 2}^2
=
{\textstyle\frac 1 4}\norm{\psi}\sb{H\sp 1}^2+A
-(B^2+{\textstyle\frac{1}{4}})\norm{\psi}^2\sb{L\sp 2}.
\end{equation}
The first inequality follows from (\ref{u-is-such}),
while the second one holds due to the bound (\ref{psi-0}).
We rewrite (\ref{q-h-0}) as the bound on $\norm{\psi}\sb{H\sp 1}^2$:
\begin{equation}
\label{q-h}
\norm{\psi}\sb{H\sp 1}^2
\le
(8B^2+2)Q(\psi)+4\mathscr{H}(\psi)-4A.
\end{equation}
When we take into account the energy and charge conservation (\ref{star}),
the inequality (\ref{q-h})
leads to the bound (\ref{a-priori-h1})
with
\begin{equation}\label{def-Lambda}
\Lambda(\phi)
=\sqrt{(8B^2+2) Q(\phi)+4\mathscr{H}(\phi)-4A}.
\end{equation}
\end{proof}

\begin{lemma}\label{lemma-u-u}
Let us assume that  
Theorem~\ref{theorem-gwp}
is true for the nonlinearities $U$ that satisfy
the following additional condition:
\begin{equation}\label{u-p-bound}
{\it For}
\quad
k=0,\,1,\,2
\quad
{\it there\ exist}
\quad
U\sb k<\infty
\quad
{\it so\ that}
\quad
\sup\sb{z\in\C}\abs{\nabla^k U(z)}\le U\sb k.
\end{equation}
Then Theorem~\ref{theorem-gwp}
is also true without this additional condition.
\end{lemma}

\begin{proof}
Fix a nonlinearity $U$
that does not necessarily satisfy (\ref{u-p-bound}).
For a particular initial data $\phi\in H\sp 1(\R)$
in Theorem~\ref{theorem-gwp},
we choose $\widetilde{U}(z)\in C\sp 2(\C)$
so that
$\widetilde{U}(z)=\widetilde{U}(\abs{z})$ for $z\in\C$
and
$\widetilde{U}(z)=U(z)$ for $\abs{z}\le\Lambda(\phi)$,
where $\Lambda(\phi)$ is defined by
(\ref{def-Lambda}).
We can choose $\widetilde{U}$
so that it satisfies (\ref{u-is-such})
with the same $A$, $B$ as $U$ does,
and also satisfies the uniform bounds
\[
\sup\sb{z\in\C}\abs{\na^k\widetilde{U}(z)}<\infty,
\qquad
k=0,1,2.
\]
By the assumption of the Lemma,
Theorem~\ref{theorem-gwp}
is true for the nonlinearity
$\widetilde F=-\nabla\widetilde{U}$
instead of $F=-\nabla U$.
Hence, there is a unique solution
$\psi(x,t)\in L\sp\infty(\R,H\sp 1)\cap C\sb b(\R\times\R)$
to the equation
\[
i\dot\psi(x,t)
=-\psi''(x,t)
-\delta(x)\widetilde{F}(\psi(0,t)),
\]
with $\psi\at{t=0}=\phi$.
By Lemma~\ref{lemma-a-priori},
$\psi$ satisfies the a priori bound (\ref{a-priori-h1})
with $\Lambda(\phi)$ defined by  (\ref{def-Lambda}).
This bound implies that
$\abs{\psi(0,t)}
\le \Lambda(\phi)$  for $t\in\R$.
Therefore, $\widetilde F(\psi(0,t))=F(\psi(0,t))$ for $t\in\R$,
and
$\psi(x,t)$
is also a solution to (\ref{snls})
with the nonlinearity $F=-\nabla U$.
\end{proof}

From now on, we shall assume in the proof of
Theorem~\ref{theorem-gwp}
that the bounds (\ref{u-p-bound}) hold true.

\begin{lemma}
\label{lemma-l8-existence}
\begin{enumerate}
\item
Let $\phi\in H^1:=H^1(\R)$.
There exists
$\tau>0$
that depends only on $U\sb 2$ in (\ref{u-p-bound})
so that there is a unique solution
$\psi\in C\sb b(\R\times[0,\tau])$
to equation (\ref{snls})
with the initial data $\psi\at{t=0}=\phi$.
\item
The map $\phi\mapsto\psi$
is continuous from $H^1$ to
$C\sb b(\R\times[0,\tau])$.
\end{enumerate}
\end{lemma}

\begin{proof}
Let us denote
the dynamical group for the free Schr\"odinger equation by
\begin{equation}\label{def-w0}
\eub{W}\sb t\phi(x)
=\frac{1}{\sqrt{2\pi t}}
\int\sb{\R} e^{i\frac{\abs{x-y}^2}{2t}}\phi(y)\,dy,
\qquad
x\in\R.
\end{equation}
For its Fourier transform, we have:
\begin{equation}
\label{fourier-transform}
\mathcal{F}\sb{x\to k}[\eub{W}\sb t\phi(x)](k)=e^{ik^2t}\hat\phi(k),
\qquad
k\in\R.
\end{equation}
Then the solution $\psi$ to (\ref{snls})
with the initial data $\psi\at{t=0}=\phi$
admits the Duhamel representation
\begin{equation}\label{duhamel}
\psi(x,t)
=\eub{W}\sb t\phi(x)
=
\eub{W}\sb t\phi(x)
+\eub{Z}\psi(x,t),
\end{equation}
where
\begin{equation}\label{def-z}
\eub{Z}\psi(x,t)
=-\int\sb{0}\sp{t}
\eub{W}\sb s\delta(x)F(\psi(0,t-s))\,ds
=
-\int\sb{0}\sp t
\frac{e^{i\frac{x^2}{2s}}}{\sqrt{2\pi s}}F(\psi(0,t-s))
\,ds.
\end{equation}
The Fourier representation (\ref{fourier-transform}) implies that
$
\eub{W}\sb t\phi(x)\in C\sb b(\R,H^1)
\subset C\sb b(\R\times\R)
$.
Further, we compute for
$\psi\sb 1$, $\psi\sb 2\in C\sb b(\R\times[0,\tau])$:
\[
\Abs{\eub{Z}\psi\sb 2(x,t)
-\eub{Z}\psi\sb 1(x,t)}
\le
\int\limits\sb{0}\sp t
\frac{
\abs{F(\psi\sb 2(0,t-s))-F(\psi\sb 1(0,t-s))}
}
{\sqrt{2\pi s}}
\,ds
\le
U\sb 2\sqrt{t}
\sup\sb{0\le s\le t} \abs{\psi\sb 2(s)-\psi\sb 1(s)},
\]
where we used (\ref{u-p-bound})
with $k=2$.
For definiteness, we set
\begin{equation}\label{def-tau-0}
\tau=\frac{1}{4U\sb 2^2}.
\end{equation}
Then the map $\psi \mapsto \eub{W}\sb{t}\phi+\eub{Z}\psi$
is contracting in the space $C\sb b(\R\times [0,\tau])$.
It follows that
equation (\ref{duhamel}) admits a unique solution
$\psi\in C\sb b(\R\times[0,\tau])$,
proving the first part of the theorem.
The second part of the theorem also follows by contraction.
\end{proof}

\section{Regularized equation}
\label{sect-regularized}
We proved that there is a unique solution
$\psi(x,t)
\in C([0,\tau]\times\R)$.
Now we are going to prove that
$\psi\in L\sp\infty(\R\sb{+},H^1)$
and moreover that
$\norm{\psi(t)}\sb{H\sp 1}$
is bounded uniformly in time.

\bigskip

Let us fix a family of functions
$\rho\sb\epsilon(x)$
approximating
the Dirac $\delta$-function.
We pick $\rho\sb 1(x)\in C\sp\infty\sb 0[-1,1]$,
nonnegative,
and such that $\int\sb{\R}\rho\sb 1(x)\,dx=1$,
and define
\begin{equation}\label{def-rho-epsilon}
\rho\sb\epsilon(x)=\frac{1}{\epsilon}\rho\sb 1\left(\frac{x}{\epsilon}\right),
\qquad
\epsilon\in(0,1),
\end{equation}
so that
\[
\supp\rho\sb\epsilon(x)\subseteq [-\epsilon,\epsilon],
\qquad\rho\sb\epsilon(x)\ge 0,
\qquad
\int\sb{\R}
\rho\sb\epsilon(x)\,dx=1.
\]
Consider the smoothed equation
with the ``mean field interaction''
\begin{equation}\label{snls-rho}
i\dot\psi(x,t)
=-\Delta\psi(x,t)-\rho\sb\epsilon(x)
F(\langle\rho\sb\epsilon,\psi(t)\rangle),
\qquad
\end{equation}
where
\[
\langle\rho\sb\epsilon,\psi(t)\rangle
=\langle\rho\sb\epsilon(\cdot),\psi(\cdot,t)\rangle
=\int\sb{\R}\rho\sb\epsilon(x)\psi(x,t)\,dx.
\]
Clearly,
equation (\ref{snls-rho})
is the Hamiltonian equation,
with the Hamilton functional
\begin{equation}
\mathscr{H}\sb\epsilon(\psi)=\int\frac{\abs{\nabla\psi}^2}{2}\,dx
+U(\langle\rho\sb\epsilon,\psi\rangle).
\end{equation}
The Hamiltonian form of (\ref{snls-rho}) is (cf.  (\ref{Ham}))
\begin{equation}\label{Hameps}
\dot\Psi\sb\epsilon=J D \mathscr{H}\sb\epsilon(\Psi\sb\epsilon).
\end{equation}
The solution $\psi\sb\epsilon$
to (\ref{snls-rho})
with the initial data $\psi\sb\epsilon\at{t=0}=\phi$
admits the Duhamel representation
\begin{equation}\label{Duep}
\psi\sb\epsilon(x,t)
=\eub{W}\sb t\phi(x)
+\eub{Z}\sb\epsilon\psi\sb\epsilon(x,t),
\end{equation}
where
\begin{equation}\label{def-z-epsilon}
\eub{Z}\sb\epsilon\psi\sb\epsilon(x,t)
=
-\int\sb{0}\sp{t}
\eub{W}\sb s
\rho\sb\epsilon(x)
F(\langle\rho\sb\epsilon,\psi\sb\epsilon(t-s)\rangle)
\,ds.
\end{equation}

\begin{lemma}[Local well-posedness]
\label{lemma-local-existence}
\begin{enumerate}
\item
For any $\epsilon\in (0,1)$,
there exists $\tau\sb\epsilon>0$
that depends on $\epsilon$ and on $U\sb 2$
from (\ref{u-p-bound})
so that there is a unique solution
$\psi\sb\epsilon\in C\sb b([0,\tau\sb\epsilon],H^1)$
to equation (\ref{snls-rho})
with $\psi\sb\epsilon\at{t=0}=\phi$.
\item
For each $t\le\tau\sb\epsilon$,
the map
$\eub{U}\sb\epsilon(t):
\phi=\psi\sb\epsilon(0)\mapsto \psi\sb\epsilon(t)$
is continuous in $H^1$.

\item
The values of the functionals
$\mathscr{H}\sb\epsilon$ and $Q$
on solutions to (\ref{snls-rho})
are conserved in time.
\end{enumerate}
\end{lemma}

\begin{proof}
\begin{enumerate}
\item
For $\psi\sb 1$, $\psi\sb 2\in C\sb b([0,\tau\sb\epsilon],H^1)$,
we compute:
\begin{eqnarray*}
&&\norm{\eub{Z}\sb\epsilon\psi\sb{2}(\cdot,t)
-\eub{Z}\sb\epsilon\psi\sb{1}(\cdot,t)}\sb{H\sp 1}
\\
&&=\norm{
\int\sb{0}\sp{t}
\eub{W}\sb s
\rho\sb\epsilon
F(\langle\rho\sb\epsilon,\psi\sb{2}(t-s)\rangle)
-F(\langle\rho\sb\epsilon,\psi\sb{1}(t-s)\rangle)
\,ds
}\sb{H\sp 1}
\\
&&\le
\int\sb{0}\sp{t}
\Norm{
\eub{W}\sb s
\rho\sb\epsilon
}\sb{H\sp 1}
\abs{
F(\langle\rho\sb\epsilon,\psi\sb{2}(t-s)\rangle)
-F(\langle\rho\sb\epsilon,\psi\sb{1}(t-s)\rangle)
}
\,ds.
\end{eqnarray*}
The first factor under the integral sign
is bounded uniformly for $0<s\le t$:
\begin{eqnarray}
\nonumber
\norm{
\eub{W}\sb s
\rho\sb\epsilon
}\sb{H\sp 1\sb x}
=\frac{1}{\sqrt{2\pi}}
\Norm{
\sqrt{1+k^2}
e^{ik^2 s/2}
\widehat{\rho\sb\epsilon}(k)
}\sb{L\sp 2\sb k}
=
\Norm{
\rho\sb\epsilon
}\sb{H\sp 1}.
\nonumber
\end{eqnarray}
Taking this into account, we get:
\begin{eqnarray*}
\norm{\eub{Z}\sb\epsilon\psi\sb{2}(\cdot,t)
-\eub{Z}\sb\epsilon\psi\sb{1}(\cdot,t)}\sb{H\sp 1}
&\le&
\norm{\rho\sb\epsilon}\sb{H\sp 1}
\int\sb{0}\sp{t}
\abs{
F(\langle\rho\sb\epsilon,\psi\sb{2}(t-s)\rangle)
-F(\langle\rho\sb\epsilon,\psi\sb{1}(t-s)\rangle)
}
\,ds
\\
&\le&
t U\sb 2\norm{\rho\sb\epsilon}\sb{H\sp 1}
\sup\sb{s\in[0,t]}
\abs{
\langle\rho\sb\epsilon,\psi\sb{2}(s)-\psi\sb{1}(s)\rangle}.
\end{eqnarray*}
Therefore, the map $\psi\mapsto\eub{W}\sb t\phi+\eub{Z}\sb\epsilon\psi$
is contracting if we choose, for definiteness,
\begin{equation}\label{def-tau-e}
\tau\sb\epsilon=\frac{1}{4 U\sb 2\norm{\rho\sb\epsilon}\sb{H\sp 1}}.
\end{equation}

\item
The continuity of the mapping
$\eub{U}\sb\epsilon(t)$
also follows from the contraction argument.

\item
It suffices to prove
the conservation of the values of
$\mathscr{H}\sb\epsilon(\psi\sb\epsilon(t))$
and $Q(\psi\sb\epsilon(t))$
for $\phi\in H^2:=H^2(\R)$ since the functionals
are continuous on $H^1$.
For $\phi\in H^2$,
the corresponding solution belongs to the space
$C\sb b([0,\tau\sb\epsilon],H^2)$
by
the Duhamel representation (\ref{Duep}).
Then
the energy and charge conservation follows
by the Hamiltonian structure (\ref{Hameps}).
Namely, the differentiation of the Hamilton functional
gives by the chain rule,
\begin{equation}\label{eneps}
\displaystyle
\frac d{dt}\mathscr{H}\sb\epsilon(\Psi\sb\epsilon(t))=
\langle
D\mathscr{H}\sb\epsilon(\Psi\sb\epsilon(t)),\dot\Psi\sb\epsilon(t)
\rangle
=
\langle
D\mathscr{H}\sb\epsilon(\Psi\sb\epsilon(t)),JD\mathscr{H}\sb\epsilon(\Psi\sb\epsilon(t))
\rangle
=0
\end{equation}
since the Fr\'echet derivative
$D\mathscr{H}\sb\epsilon(\Psi\sb\epsilon(t))
=-\Delta\Psi\sb\epsilon(\cdot,t)
-\rho\sb\epsilon(\cdot)F(\langle\rho\sb\epsilon,\Psi\sb\epsilon(t)) \rangle)
$
belongs to $L^2(\R)$ for $t\in[0,\tau\sb\epsilon]$.
Similarly, the charge conservation
follows  by the differentiation,
\begin{eqnarray}\label{qeps}
\displaystyle
\frac d{dt}Q (\Psi\sb\epsilon(t))&=&
\langle
D Q(\Psi\sb\epsilon(t)),\dot\Psi\sb\epsilon(t)
\rangle
=
\langle
D Q(\Psi\sb\epsilon(t)),J D\mathscr{H}\sb\epsilon(\Psi\sb\epsilon(t))
\rangle
\nonumber\\
&=&
\langle
\Psi\sb\epsilon(x,t),J\Delta\Psi\sb\epsilon(x,t)\rangle
-
\langle\Psi\sb\epsilon(x,t),J\rho\sb\epsilon(x)
F(\langle\rho\sb\epsilon,\Psi\sb\epsilon(t)\rangle)
\rangle.
\end{eqnarray}
Here
$\Psi\sb\epsilon(x,t),J\Delta\Psi\sb\epsilon(x,t)\rangle
=\na\Psi\sb\epsilon(x,t),J\na\Psi\sb\epsilon(x,t)\rangle=0$,
and also
\begin{eqnarray}
\langle\Psi\sb\epsilon(x,t),J\rho\sb\epsilon(x)
F(\langle\rho\sb\epsilon,\Psi\sb\epsilon(t)\rangle)\rangle
&=&
\int\Psi\sb\epsilon(x,t)
\cdot[J\rho\sb\epsilon(x)F(\langle\rho\sb\epsilon,\Psi\sb\epsilon(t)\rangle)]
\,dx\nonumber\\
&=&
\langle\rho\sb\epsilon,\Psi\sb\epsilon(t)\rangle
\cdot[J F(\langle\rho\sb\epsilon,\Psi\sb\epsilon(t)\rangle)]
=0.
\end{eqnarray}
Here ``$\cdot$'' stands for the real scalar product in $\R^2$,
and $Z\cdot[J F(Z)]=0$ for $Z\in\R^2$
since
$F(Z)=a(\abs{Z})Z$ with $a(\abs{Z})\in\R$ by (\ref{def-a}).
\end{enumerate}
\end{proof}

\begin{corollary}[Global well-posedness]
\label{corollary-global}
\begin{enumerate}
\item
For any $\epsilon>0$, $\epsilon\le 1$,
there exists a unique solution
$\psi\sb\epsilon\in C(\R,H^1)$
to equation (\ref{snls-rho})
with $\psi\sb\epsilon\at{t=0}=\phi$.

The $H\sp 1$-norm of $\psi\sb\epsilon$
is bounded uniformly in time:
\begin{equation}\label{eps}
\sup\limits\sb{t\in\R}\norm{\psi\sb\epsilon(t)}\sb{H\sp 1}
\le\Lambda\sb\epsilon(\phi),
\qquad
t\in\R,
\end{equation}
where
\begin{equation}
\Lambda\sb\epsilon(\phi)
=
\sqrt{(8B^2+2)Q(\phi)
+4\mathscr{H}\sb\epsilon(\phi)-4A}.
\end{equation}
\item
For each $t\ge 0$,
the map
$\eub{U}\sb\epsilon(t):
\psi\sb\epsilon(0)\mapsto \psi\sb\epsilon(t)$
is continuous in $H^1$.
\end{enumerate}
\end{corollary}

\begin{proof}
\begin{enumerate}
\item
The existence and uniqueness of the solution
$\psi\sb\epsilon\in C\sb b([0,\tau\sb\epsilon],H^1)$ follow
from Lemma~\ref{lemma-local-existence} ({\it i}).
The bound on the value of the $H\sp 1$-norm
of $\psi\sb\epsilon(t)$
is obtained as in Lemma~\ref{lemma-a-priori}.
Namely, noting that
\[
U(\langle\rho,\psi\sb\epsilon\rangle)
\ge
A-B\langle\rho,\psi\sb\epsilon\rangle^2
\ge
A-B\sup\sb{x\in\R}\abs{\psi\sb\epsilon}^2
\ge
A
-B^2\norm{\psi}\sb{L\sp 2}^2
-\frac 1 4\norm{\psi'}\sb{L\sp 2}^2
\]
and using the energy and charge conservation
proved in Lemma~\ref{lemma-local-existence} ({\it iii}),
we conclude that
\[
(2B^2+{\textstyle\frac 1 2})Q(\phi)
+\mathscr{H}\sb\epsilon(\phi)
=
(2B^2+{\textstyle\frac 1 2})Q(\psi\sb\epsilon)
+\mathscr{H}\sb\epsilon(\psi\sb\epsilon)
\ge
A+{\textstyle\frac 1 4}
\norm{\psi\sb\epsilon}\sb{H\sp 1}^2,
\]
so that
\begin{equation}\label{M1}
\norm{\psi\sb\epsilon}\sb{H\sp 1}^2
\le
(8B^2+2)Q(\phi)
+4\mathscr{H}\sb\epsilon(\phi)-4A.
\end{equation}
By (\ref{def-tau-e}),
the time span $\tau\sb\epsilon$ depends only on
$\norm{\rho\sb\epsilon}\sb{H\sp 1}$ and $U\sb 2$.
Hence, the bound  (\ref{eps}) allows us to
extend the solution to
$t\in[\tau\sb\epsilon,2\tau\sb\epsilon]$.
The bound  (\ref{eps}) 
for $t\in[0,2\tau\sb\epsilon]$
follows from  (\ref{M1}) 
by the energy and charge conservation
proved in Lemma~\ref{lemma-local-existence} ({\it iii}).
We conclude by induction that the solution exists
and the bound  (\ref{eps}) holds
for all $t\in\R$.

\item
The continuity of the mapping
$\eub{U}\sb\epsilon(t):
\psi\sb\epsilon(0)\mapsto \psi\sb\epsilon(t)$
for all $t\ge 0$ follows from its continuity for small times
by dividing the interval $[0,t]$ into small time intervals.
\end{enumerate}
\end{proof}

\section{Convergence of regularized solutions}
\label{sect-convergence}

\begin{lemma}
\label{lemma-uniform}
Let $\tau$
and $\psi\in C\sb b(\R\times[0,\tau])$
be as in Lemma~\ref{lemma-l8-existence},
and let $\psi\sb\epsilon\in C(\R\sb{+},H^1)$
be as in Corollary~\ref{corollary-global}.
Then for any finite $R>0$
\begin{equation}\label{W}
\psi\sb\epsilon(x,t)
\mathop{\rightrightarrows}\limits\sb{\epsilon\to 0}
\psi (x,t),
\qquad
\abs{x}\le R,
\quad
0\le t\le\tau.
\end{equation}
\end{lemma}

\begin{proof}
We have
\begin{eqnarray}
\psi\sb\epsilon(x,t)
=\eub{W}\sb t\phi(x)
+\int\sb{0}\sp t
\eub{W}\sb s
\rho\sb\epsilon(x)
F(\langle\rho\sb\epsilon,\psi\sb\epsilon(t-s)\rangle)
\,ds,
\\
\psi(x,t)
=\eub{W}\sb t\phi(x)
+\int\sb{0}\sp t
\eub{W}\sb s
\delta(x)
F(\psi(0,t-s))
\,ds.
\end{eqnarray}
Taking the difference of these equations
and regrouping the terms,
we can write:
\begin{eqnarray}
\psi\sb\epsilon(x,t)-\psi(x,t)
=\int\sb{0}\sp t
\eub{W}\sb s
\rho\sb\epsilon(x)
\left(
F(\langle\rho\sb\epsilon,\psi\sb\epsilon(t-s)\rangle)
-
F(\psi(0,t-s))
\right)\,ds
\nonumber
\\
+\int\sb{0}\sp{t}
[\eub{W}\sb s
\rho\sb\epsilon(x)-\eub{W}\sb s\delta(x)
]
F(\psi(0,t-s))
\,ds.
\label{lhs}
\end{eqnarray}
Let us analyze the first term
in the right-hand side of (\ref{lhs}).
It is bounded by
\begin{eqnarray}
\Abs{\int\sb{0}\sp t
\frac{e^{i\frac{(x-y)^2}{2s}}}{\sqrt{2\pi s}}
\rho\sb\epsilon(y)\,dy\,ds}
\sup\sb{0\le s\le t}
\abs{F(\langle\rho\sb\epsilon,\psi\sb\epsilon(s)\rangle)
-
F(\psi(0,s))}
\nonumber
\\
\le
\Abs{\int\sb{0}\sp t
\frac{ds}{\sqrt{2\pi s}}
}
U\sb 2
\sup\sb{\abs{x}\le\epsilon,\,0\le s\le t}
\abs{\psi\sb\epsilon(x,s)-\psi(x,s)}
\nonumber
\\
\le
\sqrt{\frac{2t}{\pi}}U\sb 2
\sup\sb{\abs{x}\le\epsilon,\,0\le s\le t}
\abs{\psi\sb\epsilon(x,s)-\psi(x,s)}
\nonumber
\\
\le
\frac 1 2
\sup\sb{\abs{x}\le\epsilon,\,0\le s\le t}
\abs{\psi\sb\epsilon(x,s)-\psi(x,s)},
\label{smaller-than-half}
\end{eqnarray}
where in the last inequality we used (\ref{def-tau-e}).
Setting
$M\sb{R,\tau}=\sup\sb{\abs{x}\le R,\,0\le t\le\tau}
\abs{\psi\sb\epsilon(x,t)-\psi(x,t)}$,
we can rewrite (\ref{lhs}) as
\[
M\sb{R,\tau}
\le
\frac 1 2 M\sb{R,\tau}
+
\sup\sb{\abs{x}\le R,\,0\le t\le\tau}
\int\sb{0}\sp{t}
[\eub{W}\sb s
\rho\sb\epsilon(x)-\eub{W}\sb s\delta(x)
]
F(\psi(0,t-s))
\,ds.
\]
Therefore,
\begin{equation}\label{sup-delta-psi}
M\sb{R,\tau}
\le
2\sup\sb{\abs{x}\le R,\,0\le t\le\tau}
\int\sb{0}\sp{t}
\int
\frac{e^{i\frac{(x-y)^2}{2s}}}{\sqrt{2\pi s}}
[
\rho\sb\epsilon(y)-\delta(y)
]\,dy\,F(\psi(0,t-s))ds.
\end{equation}
We claim that the right-hand side tends to zero as
$\epsilon\to 0$.
To prove this, we split the integral
into two pieces:
\begin{eqnarray}
I\sb 1(\delta,\epsilon)
=\int\sb{\delta}\sp{t}\int
\frac{e^{i\frac{(x-y)^2}{2s}}}{\sqrt{2\pi s}}
[
\rho\sb\epsilon(y)-\delta(y)
]\,dy\,F(\psi(0,t-s))ds,
\label{second-guy}
\\
I\sb 2(\delta,\epsilon)
=\int\sb{0}\sp \delta\int
\frac{e^{i\frac{(x-y)^2}{2s}}}{\sqrt{2\pi s}}
[
\rho\sb\epsilon(y)-\delta(y)
]\,dy\,F(\psi(0,t-s))ds,
\label{third-guy}
\end{eqnarray}
where $\delta\in(0,t)$ is yet to be chosen.
Let us analyze the term (\ref{second-guy}):
\begin{equation}
\abs{I\sb 1(\delta,\epsilon)}
\le
C U\sb 0
\sup\sb{s\ge\delta,\abs{x}\le R}
\Abs{
\int\sb{\abs{y}<\epsilon}
\frac{e^{i\frac{(x-y)^2}{2s}}}{\sqrt{2\pi s}}
[
\rho\sb\epsilon(y)-\delta(y)
]
\,dy}.
\end{equation}
Since $s\ge\delta>0$ and $\abs{x}\le R$,
the function
$\frac{e^{i\frac{(x-y)^2}{2s}}}{\sqrt{2\pi s}}$
is Lipschitz in $y\in [-\epsilon,\epsilon]$, uniformly in all the parameters.
Therefore,
\begin{equation}
\int\sb\R
\frac{e^{i\frac{(x-y)^2}{2s}}}{\sqrt{2\pi s}}
[
\rho\sb\epsilon(y)-\delta(y)
]
\,dy
\to 0,
\qquad
\epsilon\to 0,
\end{equation}
uniformly in the parameters.
We conclude that
\begin{equation}\label{i-1}
\lim\sb{\epsilon\to 0}
I\sb 1(\delta,\epsilon)=0,
\end{equation}
for any fixed $\delta>0$.
We then bound (\ref{third-guy})
uniformly by
\[
I\sb 2(\delta,\epsilon)
\le C U\sb 0
\int(\rho\sb\epsilon(y)+\delta(y))\,dy
\int\sb{0}\sp \delta\frac{ds}{\sqrt{s}}
\le C\sqrt{\delta},
\]
with $C$ independent of $\epsilon$.
Now apparently
the right-hand side of  (\ref{sup-delta-psi})
tends to zero as $\epsilon\to 0$.
\end{proof}

\section{Well-posedness in energy space}
\label{sect-wp}

\begin{lemma}[Local well-posedness]
\label{lemma-h1}
There is a unique  solution
$\psi\in L^\infty([0,\tau],H^1(\R))\cap C\sb b(\R\times[0,\tau])$
to equation (\ref{snls}) with $\psi\at{t=0}=\phi$,
where $\tau$ is as in (\ref{def-tau-0}).
\end{lemma}

\begin{proof}
The unique solution $\psi\in C\sb b(\R\times[0,\tau])$
is constructed in Lemma~\ref{lemma-l8-existence}.
According to (\ref{eps}) and (\ref{W}),
\begin{equation}\label{mono-h}
\norm{\psi(t)}\sb{H\sp 1}
\le\liminf \sb{\epsilon\to 0}\norm{\psi\sb\epsilon(t)}\sb{H\sp 1}
\le \Lambda(\phi),
\qquad
0\le t\le\tau.
\end{equation}
\end{proof}

\begin{lemma}
The values of the functionals
$\mathscr{H}$ and $Q$
are conserved in time for $t\in[0,\tau]$.
\end{lemma}

\begin{proof}
The convergence
(\ref{W}) and the bounds (\ref{eps}) imply that
\begin{equation}\label{charge-conservation}
Q(\psi(t))
=
\frac 1 2
\norm{\psi(t)}\sb{L\sp 2}^2
\le
\frac 1 2
\lim\sb{\epsilon\to 0}\norm{\psi\sb\epsilon(t)}\sb{L\sp 2}^2
=Q(\phi),
\end{equation}
where we used the conservation of $Q$
for the approximate solutions $\psi\sb\epsilon$
(Lemma~\ref{lemma-local-existence}).
The same argument
applied to the initial data $\psi\at{t=t\sb 0}$
with any $t\sb 0\in (0,\tau)$
and combined with the uniqueness of the solution,
allows to conclude that
$Q(\psi(t))$ is monotonically non-increasing
when time changes from $0$ to $\tau$.
Instead,
solving the Schr\"odinger equation backwards in time
and using the uniqueness of solution,
we can as well conclude that
$Q(\psi(t))$ is monotonically non-decreasing
when time changes from $0$ to $\tau$.
This proves that
$Q(\psi(t))=\const$ for $t\in[0,\tau]$.

To prove the conservation of $\mathscr{H}(\psi(t))$,
we will need the relation
\begin{equation}\label{lim-u-good}
\lim\sb{\epsilon\to 0}
U(\langle\rho\sb\epsilon,\psi\sb\epsilon\rangle)
=U(\psi(0,t)).
\end{equation}
This relation
follows from continuity of the potential $U$
and from
\begin{equation}\label{lim-is-lim}
\lim\sb{\epsilon\to 0}
\langle\rho\sb\epsilon,\psi\sb\epsilon(t)\rangle
=
\lim\sb{\epsilon\to 0}
\langle\rho\sb\epsilon,(\psi\sb\epsilon(t)-\psi(t))\rangle
+
\lim\sb{\epsilon\to 0}
\langle\rho\sb\epsilon,\psi(t)\rangle
=\psi(0,t),
\end{equation}
where
$\lim\sb{\epsilon\to 0}
\langle\rho\sb\epsilon,(\psi\sb\epsilon(t)-\psi(t))\rangle
=0$
since $\psi\sb\epsilon$
approaches $\psi$ uniformly for $0\le t\le\tau$
and $\abs{x}\le R$ (including $x=0$),
while
$\lim\sb{\epsilon\to 0}
\langle\rho\sb\epsilon,\psi(t)\rangle
=\psi(0,t)
$
since $\psi$
is continuous in $x$
(due to the finiteness of $H\sp 1$-norm
of $\psi$ that follows from (\ref{mono-h})).
We have:
\[
\mathscr{H}(\psi(t))
=
\frac{\norm{\na \psi(x,t)}\sb{L^2}^2}{2}
+U(\psi(0,t))
\le\lim\sb{\epsilon\to 0}
\left\{
\frac{\norm{\na\psi\sb\epsilon(x,t)}\sb{L^2}^2}{2}
+U(\langle\rho\sb\epsilon,\psi\sb\epsilon\rangle)
\right\}
=\mathscr{H}(\phi),
\]
where
we used the relation (\ref{lim-u-good})
and (\ref{W}).
We also used the conservation of
the values of the functional $\mathscr{H}\sb\epsilon$
for the approximate solutions $\psi\sb\epsilon$
(see Lemma~\ref{lemma-local-existence}).
Proceeding just as with $Q(\psi(t))$ above,
we conclude that  $\mathscr{H}(\psi(t))=\const$
for $0\le t\le\tau$.
\end{proof}

\begin{corollary}[Global well-posedness]
\label{corollary-global-energy}
There is a unique  solution
$\psi\in L^\infty(\R,H^1(\R))\cap C\sb b(\R\times\R)$
to equation (\ref{snls}) with $\psi\at{t=0}=\phi$.
The values of the functionals
$\mathscr{H}$ and $Q$
are conserved in time.
\end{corollary}

\begin{proof}
The solution
$\psi\in L\sp\infty([0,\tau],H\sp 1)$
constructed in Lemma~\ref{lemma-h1}
exists for $0\le t\le\tau$,
where the time span $\tau$
defined in (\ref{def-tau-0}) depends only on $U\sb 2$ from (\ref{u-p-bound}).
Hence, the bound  (\ref{a-priori-h1})
at $t=\tau$
allows us to
extend the solution $\psi$ constructed in Lemma~\ref{lemma-h1}
to the time interval $[\tau,2\tau]$.
We proceed by induction.
\end{proof}

For the conclusion of Theorem~\ref{theorem-gwp},
it remains to prove that $\psi\in C(\R,H\sp 1(\R))$.
This follows from the next two lemmas.

\begin{lemma}\label{lemma-h1-weak}
$\psi\in C(\R, H^1\sb{weak}(\R))$.
\end{lemma}

\begin{proof}
Fix $f\in H^{-1}(\R)$ and pick any $\delta>0$.
Since $H\sp 1$ is dense in $H\sp{-1}$,
there exists
$g\in H\sp 1(\R)$ such that
\begin{equation}\label{qq}
\norm{f-g}\sb{H\sp{-1}}
<\frac{\delta}{4\Lambda(\phi)},
\end{equation}
where $\Lambda(\phi)$ given by (\ref{def-Lambda})
is the a priori bound on
$\norm{\psi(t)}\sb{H\sp{1}}$
proved in Lemma~\ref{lemma-a-priori}
on the grounds of the energy and the charge conservation
for $\psi(t)$.
Then
\begin{eqnarray}
&&
\abs{\langle f,\psi(t)-\psi(t\sb 0)\rangle}
\le
\abs{\langle f-g, \psi(t)-\psi(t\sb 0)\rangle}
+\abs{\langle g, \psi(t)-\psi(t\sb 0)\rangle}
\\
&&
\le
\norm{f-g}\sb{H\sp{-1}}
\big(\norm{\psi(t)}\sb{H\sp{1}}+\norm{\psi(t\sb 0)}\sb{H\sp{1}}\big)
+\norm{g}\sb{H\sp 1}
\norm{\psi(t)-\psi(t\sb 0)}\sb{H\sp{-1}}.
\label{rrr}
\end{eqnarray}
By (\ref{qq}),
the first term
in the right-hand side of (\ref{rrr})
is bounded by
$\delta/2$.
By Corollary~\ref{corollary-global-energy},
we have $\psi\in L\sp\infty(\R,H\sp 1(\R))$,
and equation (\ref{snls}) yields $\psi\in C(\R,H\sp{-1}(\R))$.
Hence, the second term
in the right-hand side of (\ref{rrr})
becomes smaller than
$\delta/2$
if $t$ is sufficiently close to $t\sb 0$.
Since $\delta>0$ was arbitrary,
this proves that
$\lim\sb{t\to t\sb 0}\langle f,\psi(t)-\psi(t\sb 0)\rangle=0$.
\end{proof}

\begin{proposition}
\label{prop-c-h-1}
$\psi\in C(\R,H\sp 1(\R))$.
\end{proposition}

\begin{proof}
Let us fix $t\sb 0\in\R$
and compute
\begin{equation}\label{diff-h1}
\lim\sb{t\to t\sb 0}\norm{\psi(t)-\psi(t\sb 0)}\sb{H\sp 1}^2
=
\lim\sb{t\to t\sb 0}
\left(
\norm{\psi(t)}\sb{H\sp 1}^2
-2\langle\psi(t),\psi(t\sb 0)\rangle\sb{H\sp 1}
+\norm{\psi(t\sb 0)}\sb{H\sp 1}^2
\right).
\end{equation}
The relation
\[
\norm{\psi(t)}\sb{H^1}^2
=
2\big(Q(\psi(t))+H(\psi(t))\big)-2U(\psi(0,t)),
\]
together with the conservation of the energy and charge
and
the continuity of $\psi(0,t)$ for $t\in\R$
(see Corollary~\ref{corollary-global-energy}),
shows that
\[
\lim\sb{t\to t\sb 0}
\norm{\psi(t)}\sb{H^1}^2
=
\norm{\psi(t\sb 0)}\sb{H^1}^2.
\]
By Lemma~\ref{lemma-h1-weak},
$\lim\sb{t\to t\sb 0}
\langle\psi(t),\psi(t\sb 0)\rangle\sb{H\sp 1}
=\langle\psi(t\sb 0),\psi(t\sb 0)\rangle\sb{H\sp 1}$.
This shows that
the right-hand side of (\ref{diff-h1}) is equal to zero.
\end{proof}

Now Theorem~\ref{theorem-gwp} is proved.

\section{H\"older regularity of solution}
\label{sect-holder}

In this section, we prove Theorem~\ref{theorem-holder}.

\begin{lemma}
\label{lemma-w0}
If $\phi\in H^1$, then
$\eub{W}\sb{(\cdot)}\phi(x)\in C\sp{(1/4)}[0,\tau]$, uniformly in $x\in\R$.
\end{lemma}

\begin{proof}
Let $t$, $t'\in[0,\tau]$.
We have by the Cauchy-Schwarz inequality:
\begin{equations}
\abs{
\eub{W}\sb{t'}\phi(x)
-\eub{W}\sb t\phi(x)
}
\le
C\Abs{
\int
e^{-i k x}\Big(e^{i\frac{t' k^2}{2}}
-e^{i\frac{tk^2}{2}}\Big)
\hat\phi(k)\,dk
}
\nonumber
\\
\le
C\int
\min(1,\abs{t'-t}k^2)
\abs{\hat\phi(k)}\,dk
\le
C\left[
\int\sb\R
\frac{\min(1,\abs{t'-t}k^2)^{2}}{1+k^2}
\,dk
\right]^{\frac 1 2}
\norm{\phi}\sb{H\sp 1}.
\nonumber
\end{equations}
We bound the last integral
as follows:
\[
\int\sb\R
\frac{\min(1,\abs{t'-t}k^2)^{2}}{1+k^2}
\,dk
\le
\int\limits\sb{\abs{k}<\abs{t'-t}^{-\frac 1 2}}
\frac{\abs{t'-t}^2 k^4}{1+k^2}
\,dk
+
\int\limits\sb{\abs{k}>\abs{t'-t}^{-\frac 1 2}}
\frac{dk}{1+k^2}
\le\const\abs{t'-t}^{\frac 1 2}.
\]
\end{proof}

\begin{lemma}[Regularity of $\psi(0,t)$]
\label{lemma-t}
The unique solution $\psi\in C\sb b(\R\times[0,\tau])$
to equation (\ref{snls})
with the initial data
$\psi\at{t=0}=\phi$
constructed in Lemma~\ref{lemma-l8-existence}
satisfies
\[
\psi(0,\cdot)\in C\sp{(1/4)}[0,\tau].
\]
\end{lemma}

\begin{proof}
Due to Lemma~\ref{lemma-w0},
it suffices to consider the regularity of
$\eub{Z}\psi(0,t)$.
For
any $t,\,t'\in[0,\tau]$,
$t'<t$,
we have:
\begin{equation}\label{z-z-0}
\eub{Z}\psi(0,t')-\eub{Z}\psi(0,t)
=
\int\sb{0}\sp t
\left[
\frac{F(\psi(0,s))}{\sqrt{2\pi(t'-s)}}
-\frac{F(\psi(0,s))}{\sqrt{2\pi(t-s)}}
\right]
\,ds
+\int\sb{t}\sp{t'}
\frac{F(\psi(0,s))}{\sqrt{2\pi(t'-s)}}\,ds.
\end{equation}
The first integral in the right-hand side of (\ref{z-z-0})
is bounded by
\[
C\sb 1\int\sb{0}\sp t
\Abs{
\frac{1}{\sqrt{t'-s}}
-\frac{1}{\sqrt{t-s}}
}\,ds
\le
C\sb 2\abs{t'-t}\sp{1/2}.
\]
The second integral in the right-hand side of (\ref{z-z-0})
is also bounded by $C\abs{t'-t}\sp{1/2}$.
\end{proof}

\begin{lemma}\label{lemma-holder}
$\psi(x,\cdot)\in C\sp{(1/4)}(\R)$,
uniformly in $x\in\R$.
\end{lemma}

\begin{proof}
We have the relation
\begin{equation}\label{t-zero}
\psi(x,t)
=\eub{W}\sb{t-t\sb 0}\psi(x,t\sb 0)
+
\int\sb{0}\sp{t-t\sb 0}
\frac{e^{i\frac{x^2}{2s}}}{\sqrt{2\pi s}}F(\psi(0,t-s))\,ds.
\end{equation}
By Lemma~\ref{lemma-w0},
the first term in the right-hand side of (\ref{t-zero}),
considered as a function of time,
belongs to
$C\sp{(1/4)}(\R)$
(uniformly in $x\in\R$).
The second term
in the right-hand side of (\ref{t-zero})
is bounded by $\const\abs{t-t\sb 0}\sp{1/2}$.
This proves that
$\psi(x,\cdot)\in C\sp{(1/4)}(\R)$,
uniformly in $x$.
\end{proof}

It remains to mention that
the H\"older continuity in $x$
follows from the inclusion
$H\sp 1(\R)\subset C\sp{(1/4)}(\R)$.
Theorem~\ref{theorem-holder} is proved.

\section*{Acknowledgments}
The authors are grateful to David Stuart
for helpful remarks.

\bibliographystyle{amsalpha}
\bibliography{shd-local,all,berestycki,books,bourgain,cazenave,christ,comech,cuccagna,delort,esteban,grillakis,hormander,komech,melrose,pego,pelinovsky,phong,roudenko,seeger,smith,shatah,soffer,sogge,stein,strauss,sugim
oto,symes,tataru,uhlmann,weinstein,zworski}

\end{document}